\pdfoutput=1
\documentclass{ifacconf}[letterpaper]

\makeatletter
\let\old@ssect\@ssect  %
\makeatother

\newcommand{\pathtoparts}{parts_arXiv}
\newcommand{\pathtofigs}{figures_arXiv}
\usepackage{graphicx}  %
\usepackage{natbib}  %
\input{\pathtoparts/preamble.tex}

\makeatletter
\def\@ssect#1#2#3#4#5#6{
  \NR@gettitle{#6}  %
  \old@ssect{#1}{#2}{#3}{#4}{#5}{#6}  %
}
\makeatother
\begin{document}
\begin{frontmatter}

  \title{SPOCK: A proximal method for multistage risk-averse optimal control problems\thanksref{footnoteinfo}}

  \thanks[footnoteinfo]{The first two authors contributed equally to this paper. The work of R. Moran and P. Sopasakis has received support by EPSRC under the Queen's Collaborative Award ``BotDozer: GPU-accelerated model predictive control for autonomous heavy equipment'', which is co-funded by EquipmentShare.
    The work of A. Bodard, M. Schuurmans and P. Patrinos has received support by:
    the Research Foundation Flanders (FWO) research projects G081222N, G033822N, G0A0920N;
    EU's Horizon 2020 research and innovation programme under the Marie Skłodowska-Curie grant agreement No. 953348;
    Research Council KU Leuven C1 project No. C14/18/068;
    Fonds de la Recherche Scientifique - FNRS and the Fonds Wetenschappelijk Onderzoek - Vlaanderen under EOS project no 30468160 (SeLMA).}

  \author[KUL]{Alexander Bodard}
  \author[QUB]{Ruairi Moran}
  \author[KUL]{Mathijs Schuurmans}
  \author[KUL]{Panagiotis Patrinos}
  \author[QUB]{Pantelis Sopasakis}

  \address[QUB]{Queen's University Belfast, School of Electronics, Electrical Engineering and Computer Science, BT9 5AH, Northern Ireland, UK. (emails: \{\textbf{rmoran05}, \textbf{p.sopasakis}\}\textbf{@qub.ac.uk})}%
  \address[KUL]{KU Leuven, Department of Electrical Engineering, Kasteelpark Arenberg 10, 3001 Leuven, Belgium. (emails: \{\textbf{alexander.bodard}, \textbf{mathijs.schuurmans}, \textbf{panos.patrinos}\}\textbf{@esat.kuleuven.be})}

  \begin{abstract} %
    Risk-averse optimal control problems have gained a lot of attention in the last decade, mostly due to their attractive mathematical properties and practical importance.
They can be seen as an interpolation between stochastic and robust optimal control approaches, allowing the designer to trade-off performance for robustness and vice-versa.
Due to their stochastic nature, risk-averse problems are of a very large scale, involving millions of decision variables, which poses a challenge in terms of efficient computation.
In this work, we propose a splitting for general risk-averse problems and show how to efficiently compute iterates on a GPU-enabled hardware.
Moreover, we propose \SPOCK{} --- a new algorithm that utilizes the proposed splitting and takes advantage of the SuperMann scheme combined with fast directions from Anderson's acceleration method for enhanced convergence speed.
We implement \SPOCK{} in \texttt{Julia} as an open-source solver, which is amenable to warm-starting and massive parallelization.

  \end{abstract}


\end{frontmatter}

\section{Background and motivation}
There are two main ways to formulate optimal control problems (OCPs) for uncertain systems: (i) the \textit{stochastic/risk-neutral}~\citep{cinquemani_stochastic_estimation} approach where we assume, naively, that we have perfect information about the distributions of the involved random disturbances, and (ii) the \textit{robust/minimax} approach where we ignore any statistical information that is usually available~\cite[Ch.~3]{rawlings_mpc_book}.
Another way is the risk-averse formulation, which takes into account \textit{inexact} and \textit{data-driven} statistical information \citep{wang_risk-problems-survey} using risk measures.
A risk measure quantifies the magnitude of the right tail of a random cost, allowing the designer to choose a balance between the risk-neutral and the worst-case scenario.
We restrict our interest to \textit{coherent} risk measures \citep{shapiro_lectures-stochastic} because of the desirable mathematical properties they possess.

The challenge that multistage risk-averse optimal control problems (RAOCPs) bring, is that their cost functions involve the composition of several nonsmooth risk measures.

Risk measures have been widely used in practice and are gradually becoming popular, e.g., in reinforcement learning \citep{chow_risk-measures2017, chow_risk-measures2015}, in MPC \citep{sopasakis_risk-averse-mpc-preprint} for microgrids \citep{hans_risk-measures}, collision avoidance \citep{dixit_risk-measures, schuurmans_risk-measures}, energy management systems \citep{maree_risk-measures} and air-ground rendezvous \citep{barsi_risk-measures}.

RAOCPs are typically solved using stochastic dual dynamic programming \citep{pereira_sddp-original, shapiro_analysis-sddp} approaches such as \citep{dacosta_sddp, shapiro_ra-and-rn-sddp}, however, these turn out to be slow.
Authors in \citep{sopasakis_ra-rc-mpc} proposed an approach that allows us to use more efficient out-of-the-box optimization methods and software such as \GUROBI{} \citep{gurobi_manual} and \MOSEK{} \citep{andersen_mosek}.
These solvers use interior-point methods \citep{nemirovski_interior-point}, which usually do not scale well with the size of the problem and are hard to warm-start \citep{yildirim_warm-start-interior-point}.
However, risk-averse problems possess a rich structure that we can exploit to devise very efficient and massively parallelizable methods that can run on lockstep architectures such as GPUs~\citep{oh_lockstep_parallelisation}.

We propose a look to first-order methods because they have simple iterates that are amenable to parallelization.
In addition, the nonsmooth nature of the cost functions of risk-averse problems prohibits the use of any method that requires gradients.
The Chambolle-Pock method \citep{chambolle_method} is chosen because it is gradient-free, it scales well with the problem size, it is easy to warm-start, it has a provably decent $\mathcal{O}(1/k)$ convergence rate which is similar to the celebrated alternating direction method of multipliers (ADMM) \citep{parikh_proximal}, and it is simpler than ADMM as it requires only two invocations of the problem's linear operator instead of four. Furthermore, we can accelerate the Chambolle-Pock method by combining it with SuperMann \citep{themelis_supermann}, a Newton-type algorithm for finding fixed points of nonexpansive operators. We propose the resulting algorithm, SuperMann with Chambolle-Pock (\SPOCK).

Due to these advantages, proximal algorithms have caught the attention of researchers in the area of control \citep{stella_simple-nonlinear-mpc-panoc, stathopoulos_operator-splitting-methods, donoghue_splitting-method}.
Still, GPU computation is mostly underused in optimal control problems.
Recent successful attempts in stochastic control are \citep{chowdhury_gpu_stochastic_algo, sampathirao_stochastic_gpu}.

The contributions of this paper are, (i) we propose a problem splitting that allows us to apply the Chambolle Pock method and SuperMann, and (ii) an implementation of a new solver in \texttt{Julia}, and simulations that demonstrate the effectiveness of the method. The novel open-source solver \SPOCK{} is amenable to warm-starting and massive parallelization.

\section{Notation}
Let $\N_{[k_1, k_2]}$ denote the integers in $[k_1, k_2]$.
Let $1_n$ and $0_n$ be understood as $n$-dimensional column vectors of ones and zeroes respectively.
Let $I_n$ be understood as the $n$-dimensional identity matrix.
For $z\in\R^n$ let $[z]_+ = \max\{0,z\}$, where the max is taken element-wise.
We denote the transpose of a matrix $A$ by $A\Top$.
The adjoint of a linear operator $L:\R^n\rightarrow\R^m$ is the operator $L^*:\R^m\rightarrow\R^n$ that satisfies $y^\intercal Lx = x^\intercal L^*y$ for all $x\in\R^n,y\in\R^m$.
The dual cone $\K^*$ of a closed convex cone $\K \subseteq \R^n$ is the set $\K^* = \{ y\in\R^n | y\Top x \geq 0, \forall x\in\K \}$.
The relative interior of $\K$ is denoted by $\ri\;\K$.
We denote a second order cone of dimension $d$ by $\SOC_d$ and a translated second order cone by $\SOC_d + a$ where $a\in\R^d$.
A function $f:\R^n\rightarrow\R$ is called \textit{lower semicontinuous} (lsc) if its sublevel sets, $\{ x | f(x) \leq \alpha \}$, are closed.
The indicator of a convex set $\K$ is defined as
$
    \delta_{\K}(x) = 0 \text{, if } x\in\K \text{, and }
    \delta_{\K}(x) = \infty \text{, otherwise. }
$
The set of fixed points of an operator $T$ is denoted $\fix T$, where a fixed point $v$ is defined by $v = T(v)$.
Define $\overline{\R} = \R \cup \{\infty\}$.
The domain of a function $f:\R^n \to \overline{\R}$ is $\textbf{dom} f = \{ x: f(x) < \infty \}$.
The subdifferential of a convex function $f$ is $\partial f(x) = \{u\in\R^n \:|\: \forall y\in\R^n, (y-x)\Top u + f(x) \leq f(y)\}$.
The proximal operator of a proper lsc convex function $f: \R^n \to \overline{\R}$ with parameter $\alpha>0$ is
\begin{equation}
    \prox_{\alpha f}(x) = \argmin_v \left\{ f(v) + \tfrac{1}{2\alpha} \| v - x \|^2_2 \right\}.
\end{equation}
Note that superscript $i$ denotes scenario tree nodes, superscript $(k)$ denotes algorithm iterations, and subscript $t$ denotes time steps.

\section{The Chambolle-Pock method}\label{sec:cp-method}
The Chambolle-Pock (CP) method is a proximal method that can be used to solve optimization problems of the following form \cite[p.32]{ryu_monotone-primer-chambolle-pock}
\begin{equation}\label{eq:main-cp}
    \mathbb{P}:\minimise_{z\in\R^{n_z}} f(z) + g(L z),
\end{equation}
where $L: \R^{n_z} \to \R^{n_{\eta}}$ is a linear operator, and $f$ and $g$ are proper closed convex functions on $\R^{n_z}$ and $\R^{n_{\eta}}$, respectively. We assume strong duality, that is, $\ri\,\textbf{dom} \: f \cap \ri\,\textbf{dom} \: g(L) \neq \emptyset$, throughout this paper.
The primal-dual optimality conditions of \eqref{eq:main-cp} consist in determining a pair of $z\in\R^{n_z}$, $\eta\in\R^{n_{\eta}}$ such that
\begin{subequations}\label{eq:primal-dual-point}
    \begin{align}
        0   & \in \partial f(z) + L^* \eta,
        \label{eq:primal-dual-point:a}
        \\
        L z & \in \partial g^*(\eta).
        \label{eq:primal-dual-point:b}
    \end{align}
\end{subequations}
Our objective is to determine an approximate solution $(z, \eta)$ such that
\begin{subequations}\label{eq:primal-dual-point-approx}
    \begin{align}
        \xi_p       & \in \partial f(z) + L^* \eta,
        \label{eq:primal-dual-point-approx:a}
        \\
        L z + \xi_d & \in \partial g^*(\eta),
        \label{eq:primal-dual-point-approx:b}
    \end{align}
\end{subequations}
for some $\xi_p\in\R^{n_{z}}$ and $\xi_d\in\R^{n_{\eta}}$ of sufficiently small norm.

The CP method \cite[]{chambolle_method} recursively applies the firmly nonexpansive (FNE) operator $T$, where
$
    (
    z^{(k+1)}, \eta^{(k+1)}
    )
    =
    T(
    z^{(k)}, \eta^{(k)}
    ),
$
that is,
\begin{equation}\label{eq:chambolle-pock-operator}
    \begin{bmatrix}[l]
        z^{(k+1)}
        \\
        \eta^{(k+1)}
    \end{bmatrix}
    =
    \underbracket[0.5pt]{
        \begin{bmatrix}[l]
            \prox_{\alpha f} ( z^{(k)} {-} \alpha L^* \eta^{(k)} )
            \\
            \prox_{\alpha g^*} ( \eta^{(k)} {+} \alpha L (2 z^{(k+1)} {-} z^{(k)}) )
        \end{bmatrix}
    }_{T(z^{(k)}, \eta^{(k)})},
\end{equation}
which is a generalized proximal point method with preconditioning operator
\begin{equation}
    M(z, \eta) {}:={}
    \begin{bmatrix}
        I         & -\alpha L^*
        \\
        -\alpha L & I
    \end{bmatrix}
    \begin{bmatrix}
        z
        \\
        \eta
    \end{bmatrix}.
\end{equation}
The CP method converges if a solution exists and $0 < \alpha \|L\| < 1$, where $\|L\|$ is the operator norm $\|L\| = \max \left\{ \|L z\| : z\in Z, \|z\| \leq 1 \right\}$.

We define the residual operator $R {}\coloneqq{} \id - T$, where $\id$ is the identity operator.
We may partition the residual into its primal and dual parts as $R^{(k)} = (z^{(k)}, \eta^{(k)}) - T(z^{(k)}, \eta^{(k)}) = (r_{z}^{(k)}, r_{\eta}^{(k)})$.
We can show that the quantities
\(
\xi_p^{(k)} {}={} \tfrac{1}{\alpha}r_{z}^{(k)} - L^* r_{\eta}^{(k)},
\)
and
\(
\xi_d^{(k)} {}={} \tfrac{1}{\alpha}r_{\eta}^{(k)} - L r_{z}^{(k)}
\)
satisfy \eqref{eq:primal-dual-point-approx}. We define
\(
\xi_{(k)} = \max(\|\xi_p^{(k)}\|_{\infty}, \|\xi_d^{(k)}\|_{\infty})
\)
and terminate when $\xi_{(k)} \leq \max(\epsilon_{\text{abs}}, \epsilon_{\text{rel}} \xi_{(0)})$, for some absolute and relative tolerances $\epsilon_{\text{abs}}$ and $\epsilon_{\text{rel}}$, respectively. This termination criterion is akin to that used in \citep{sopasakis_superscs}.

The CP oracle consists of $\prox_{\alpha f}$, $\prox_{\alpha g^*}$, $L$ and $L^*$. The CP algorithm has similar convergence properties to the celebrated ADMM, but requires only two invocations of $L$ and $L^*$ instead of four.

\section{Multistage nested risk-averse OCPs}
In this section we introduce the risk-averse optimal control problem on a scenario tree.
We consider the following discrete-time linear dynamical system
\begin{equation}\label{eq:disturbed-dt-linear-system}
    x_{t+1} = A(w_t)x_t + B(w_t)u_t,
\end{equation}
for $t\in\N_{[0,N-1]}$, with state variable $x_t\in\R^{n_x}$ and control input $u_t\in\R^{n_u}$. The random disturbance $w_t\in\{1,\ldots,W\}$ is finite valued; for example, it can be an iid process or a Markov chain.

\subsection{Scenario trees}\label{sec:scenario-trees}
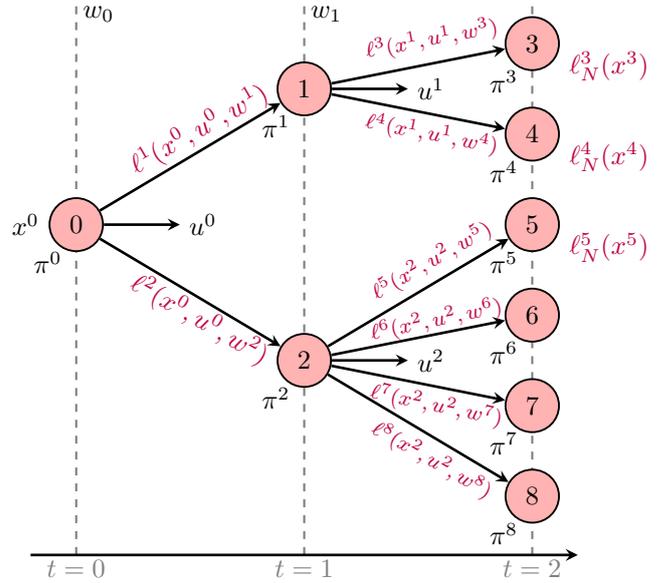
\begin{figure}
    \centering
    \begin{tikzpicture}[scale=0.6]
        \draw[thick, dashed, gray] (0,-7.3) -- (0,5);
        \draw[thick, dashed, gray] (5,-7.3) -- (5,5);
        \draw[thick, dashed, gray] (10,-7.3) -- (10,5);
        \node (node0) at (0,0)   [Nodes] {0};
        \node (node1) at (5,3)   [Nodes] {1};
        \node (node2) at (5,-3)  [Nodes] {2};
        \node (node3) at (10,4)  [Nodes] {3};
        \node (node4) at (10,2)  [Nodes] {4};
        \node (node5) at (10,0)  [Nodes] {5};
        \node (node6) at (10,-2) [Nodes] {6};
        \node (node7) at (10,-4) [Nodes] {7};
        \node (node8) at (10,-6) [Nodes] {8};
        \draw[arrow] (node0) to (node1);
        \draw[arrow] (node0) to (node2);
        \draw[arrow] (node1) to (node3);
        \draw[arrow] (node1) to (node4);
        \draw[arrow] (node2) to (node5);
        \draw[arrow] (node2) to (node6);
        \draw[arrow] (node2) to (node7);
        \draw[arrow] (node2) to (node8);
        \draw[axis] (-1,-7.3) -- (11,-7.3);
        \node (t0) at (0,-7.6) {\color{gray}$t=0$};
        \node (t1) at (5,-7.6) {\color{gray}$t=1$};
        \node (t2) at (10,-7.6) {\color{gray}$t=2$};
        \node [left=0 of node0] (x^0) {$x^0$};
        \node [right=1cm of node0] (u^0) {$u^0$};
        \draw[arrow] (node0) to (u^0);
        \node [right=1cm of node1] (u^1) {$u^1$};
        \draw[arrow] (node1) to (u^1);
        \node [right=1cm of node2] (u^2) {$u^2$};
        \draw[arrow] (node2) to (u^2);
        \node [yshift=2.8cm,right=-0.4cm of node0] (w_0) {$w_{0}$};
        \node [yshift=1.0cm,right=-0.4cm of node1] (w_1) {$w_{1}$};
        \node [rotate=30]  (ellx^0u^0w_0) at (2.7,2.1)  {\color{purple}$\ell^1(x^0,u^0,w^1)$};
        \node [rotate=-30] (ellx^0u^0w_0) at (2.7,-2.1) {\color{purple}$\ell^2(x^0,u^0,w^2)$};
        \node [rotate=11]  (ellx^1u^1w_1) at (7.8,4.1)  {\small\color{purple}$\ell^3(x^1,u^1,w^3)$};
        \node [rotate=-11] (ellx^1u^1w_1) at (7.8,2)  {\small\color{purple}$\ell^4(x^1,u^1,w^4)$};
        \node [rotate=30]  (ellx^2u^2w_1) at (7.8,-0.8) {\small\color{purple}$\ell^5(x^2,u^2,w^5)$};
        \node [rotate=12]  (ellx^2u^2w_1) at (7.9,-2.1) {\small\color{purple}$\ell^6(x^2,u^2,w^6)$};
        \node [rotate=-12] (ellx^2u^2w_1) at (7.9,-4.0) {\small\color{purple}$\ell^7(x^2,u^2,w^7)$};
        \node [rotate=-30] (ellx^2u^2w_1) at (7.8,-5.2) {\small\color{purple}$\ell^8(x^2,u^2,w^8)$};
        \node [yshift=-0.3cm,right=0 of node3] (ellx^3) {\color{purple}$\ell^3_N(x^3)$};
        \node [yshift=-0.3cm,right=0 of node4] (ellx^4) {\color{purple}$\ell^4_N(x^4)$};
        \node [yshift=-0.3cm,right=0 of node5] (ellx^5) {\color{purple}$\ell^5_N(x^5)$};
        \node [yshift=-0.5cm, left=-0.3 of node0] (pi^0) {$\pi^0$};
        \node [yshift=-0.5cm, left=-0.3 of node1] (pi^1) {$\pi^1$};
        \node [yshift=-0.5cm, left=-0.3 of node2] (pi^2) {$\pi^2$};
        \node [yshift=-0.5cm, left=-0.3 of node3] (pi^3) {$\pi^3$};
        \node [yshift=-0.5cm, left=-0.3 of node4] (pi^4) {$\pi^4$};
        \node [yshift=-0.5cm, left=-0.3 of node5] (pi^5) {$\pi^5$};
        \node [yshift=-0.5cm, left=-0.3 of node6] (pi^6) {$\pi^6$};
        \node [yshift=-0.5cm, left=-0.3 of node7] (pi^7) {$\pi^7$};
        \node [yshift=-0.5cm, left=-0.3 of node8] (pi^8) {$\pi^8$};
    \end{tikzpicture}
    \caption{Discrete RAOCP described by a scenario tree.}
    \label{fig:scenario_tree}
\end{figure}

A scenario tree is a representation of the dynamics of the system in Equation \eqref{eq:disturbed-dt-linear-system} over a finite number of stages given that the control actions are determined based on the system states in a causal fashion \citep{hoyland_scenario-trees-intro, dupavcova_scenario-trees-intro}, as depicted in Figure~\ref{fig:scenario_tree}. Each time period is called a \textit{stage}, $t$, and the number of time periods is called the \textit{horizon}, $N$. Each possible realization of Equation \eqref{eq:disturbed-dt-linear-system} at each stage is called a \textit{node}.
We enumerate the nodes of a tree with $i$, where $i=0$ is the \textit{root} node which corresponds to the initial state of the system. The nodes at subsequent stages for $t\in\N_{[0, N]}$ are denoted by $\nodes(t)$. Let $\nodes(t_1, t_2) = \cup_{t=t_1}^{t_2}\nodes(t)$, where $0\leq t_1\leq t_2\leq N$. The unique \textit{ancestor} of a node $i\in\nodes(1, N)$ is denoted by $\anc(i)$ and the set of \textit{children} of $i\in\nodes(t)$ for $t\in\N_{[0, N-1]}$ is $\ch(i) \subseteq \nodes(t+1)$. Each node is associated with a probability $\pi^i>0$ of occurring, where $\sum_{i\in\nodes(t)} \pi^i = 1$ for $t\in\N_{[0,N]}$.

\subsubsection{Dynamics.}
The system dynamics on the scenario tree is described by
\begin{equation}\label{eq:finite-horizon-evolution-2}
    x^{i_+} = A^{i_+}x^i + B^{i_+}u^i,
\end{equation}
where $i\in\nodes(0, N-1)$ and $i_+\in\ch(i)$. The state $x^0$ at the root node is assumed to be known.

\subsubsection{Constraints.}
Here we assume that the states and inputs must satisfy the convex constraints
\begin{subequations}\label{eq:constraints}
    \begin{align}
        \Gamma_x^i x^i + \Gamma_u^i u^i \in{} & C^i,
        \\
        \Gamma_N^j x^j \in{}                  & C_N^j,
    \end{align}
\end{subequations}
for $i\in\nodes(0, N-1)$ and $j\in\nodes(N)$, where $C^i$ and $C_N^j$ are closed convex sets, $\Gamma_x^i$ and $\Gamma_u^i$ are the state-input constraint matrices and $\Gamma_N^j$ is the terminal state constraint matrix.

\subsubsection{Quadratic costs.}
In discrete-time, discrete sample space stochastic OCPs, there is a cost associated with each node of the scenario tree.
The costs that are used to build optimal control problems in this paper are quadratic functions. On the scenario tree, for $i\in\nodes(0,N-1)$ and $i_+\in\ch(i)$, let $\ell^{i_+}(x^i, u^i, w^{i_+}) = x\iTop Q^{i_+} x^i + u\iTop R^{i_+} u^i$ be the stage cost function. For $j\in\nodes(N)$, the terminal cost is $\ell_N^j(x^j) = x\jTop Q^j_N x^j$.

\subsection{Risk measures}\label{sec:risk-measures}
Given a discrete sample space $\Omega=\{1,\ldots, n\}$  and a probability vector $\pi\in\R^n$, a random cost $Z:\Omega\to\R$ can be identified by a vector $Z=(Z^1, \ldots, Z^n) \in \R^n$.
The expectation of $Z$ with respect to $\pi$ is simply $\E^{\pi}[Z] {}={} \pi\Top Z$, and the maximum of $Z$ is defined as $\max Z = \max\{Z^i, i\in\N_{[1, n]}\}$.

A risk measure can be used to reflect the inexact knowledge of the probability measure.
A risk measure is an operator, $\rho:\R^n\to\R$, that maps a random cost $Z\in\R^n$ to a characteristic index $\rho(Z)$ that quantifies the magnitude of its right tail. A risk measure is said to be \textit{coherent} \cite[Sec. 6.3]{shapiro_lectures-stochastic} if (i) it is convex, (ii) it is monotone, that is, $\rho(Z)\leq\rho(Z')$ whenever $Z\leq Z'$ element-wise, (iii) $\rho$ is translation equivariant, that is, $\rho(Z+c1_n)=\rho(Z) + c$ for all $Z\in\R^n$ and $c\in\R$, (iv) $\rho$ is positively homogeneous. Coherent risk measures admit the following dual representation
\begin{equation}\label{eq:risk-dual-repr}
    \rho(Z) = \max_{\mu\in\mathcal{A}}\E^{\mu}[Z],
\end{equation}
where the set $\mathcal{A}$ is a closed convex set known as the \textit{ambiguity set} of $\rho$. An ambiguity set is a quantifiable way of indicating the uncertainty about the true probability distribution of $\pi$.
Equation \eqref{eq:risk-dual-repr} suggests that coherent risk measures can be seen as the worst-case expectations of $Z$, with respect to a probability vector $\mu$ that is taken from $\mathcal{A}$. Clearly, $\E^\pi$ is a coherent risk measure with $\mathcal{A}=\{\pi\}$ and $\max$ is also coherent with $\mathcal{A}$ being the probability simplex.

A popular coherent risk measure is the average value-at-risk with parameter $a\in[0,1]$, denoted by $\AVaR_a$. It is defined as
\begin{equation}
    \AVaR_a[Z] =
    \begin{cases}
        \min\limits_{t\in \R^n}
        \left\{
        t + \tfrac{1}{a}\mathbb{E}^\pi[Z-t]_{+}
        \right\}, & a\neq 0,
        \\
        \max[Z],  & a=0.
    \end{cases}
\end{equation}
The ambiguity set of $\AVaR_a$ is
\begin{equation}
    \mathcal{A}_a^{\rm avar}(\pi) = \left\{\mu \in \R^n \left| \sum_{i=1}^n \mu^i = 1, 0 \leq \mu^i \leq \tfrac{\pi^i}{a} \right.\right\}.
\end{equation}

\

\subsection{Conic representation of risk measures}
Coherent risk measures also admit the following \textit{conic representation}
\begin{equation}\label{eq:primal_risk}
    \rho[Z] = \max_{\mu \in \R^n, \nu \in \R^{n_{\nu}}}\
    \left\{ \mu\Top Z \:|\: b - E\mu - F\nu \in \K \right\},
\end{equation}
where $\K$ is a closed convex cone.
Provided strong duality holds (i.e., if there exist $\mu$ and $\nu$ so $b - E\mu - F\nu \in \ri \K$), the dual form of the risk measure in Equation \eqref{eq:primal_risk} is
\begin{equation}\label{eq:dual_risk}
    \rho[Z] = \min_{y}\ \left\{ y\Top b {}\mid{} E\Top y = Z, F\Top y = 0, y \in \K^* \right\}.
\end{equation}
For $\AVaR_{a}$, $n_{\nu}=0$,
$b =
    \begin{bmatrix}
        \pi\Top & 0_n\Top & 1
    \end{bmatrix}\Top,$
$E =
    \begin{bmatrix}
        a I_n & -I_n & 1_n
    \end{bmatrix}\Top,$
and
$\K = \R^{2n}_{+} \times \{0\}.$

\subsection{Conditional risk mappings}
Every node $i\in\nodes(t+1), t\in\N_{[0,N-1]}$ is associated with a cost value $Z^i = \ell^i(x^{\anc(i)}, u^{\anc(i)}, w^i)$. For each $t\in\N_{[0,N-1]}$ we define a random variable $Z_t := (Z^i)_{i\in\nodes(t+1)}$ on the probability space $\nodes(t+1)$.
Every node $j\in\nodes(N)$ is associated with a cost value $Z^j = \ell^j(x^j)$. We define the random variable $Z_N := (Z^j)_{j\in\nodes(N)}$ on the probability space $\nodes(N)$.

We define $Z^{[i]} := (Z^{i_+})_{i_+\in\ch(i)}, i\in\nodes(0,N-1)$. This partitions the variable $Z_t = (Z^{[i]})_{i\in\nodes(t)}$ into groups of nodes which share a common ancestor.

Let $\rho: \R^{|\ch(i)|} \to \R$ be risk measures on the probability space $\ch(i)$. For every stage $t \in \N_{[0,N-1]}$ we may define a conditional risk mapping at stage $t$, $\rho_{|t}: \R^{|\nodes(t+1)|} \to \R^{|\nodes(t)|}$, as follows \cite[Sec. 6.8.2]{shapiro_lectures-stochastic}
\begin{equation}
    \rho_{|t}[Z_t] {}:={} \left(\rho^i[Z^{[i]}]\right)_{i\in\nodes(t)}.
\end{equation}

\section{Proposed Splitting}
\subsection{Original problem}
A risk-averse optimal control problem with horizon $N$ is defined via the following multistage nested formulation \cite[Sec. 6.8.1]{shapiro_lectures-stochastic}
\begin{multline}\label{eq:original-problem}
    V^{\star} {}:={} \inf_{u_0} \rho_{|0} \bigg[ Z_0
        + \inf_{u_1} \rho_{|1} \Big[ Z_1 \\
            + \cdots
            + \inf_{u_{N-1}} \rho_{|N-1} \big[ Z_{N-1}
                + Z_{N} \big] \cdots \Big]  \bigg],
\end{multline}
subject to \eqref{eq:finite-horizon-evolution-2} and \eqref{eq:constraints}.
Note that the infima in \eqref{eq:original-problem} are taken element-wise and with respect to the control laws $u_t := (u^i)_{i\in\nodes(t)}$ for $t\in\N_{[0, N-1]}$.
This problem is decomposed using \textit{risk-infimum interchangeability} and \textit{epigraphical relaxation} in \cite[]{sopasakis_ra-rc-mpc}, allowing us to cast the original problem as the following minimization problem with conic constraints
\begin{subequations}\label{eq:min-s0-v1}
    \begin{equation}
        \hspace*{-5.6cm}
        \minimise_{\substack{
                (x^{i})_{i}, (x^{j})_{j}, (u^{i})_{i}, (y^{i})_{i},
                \\
                (\tau^{i_+})_{i_+}, (s^{i})_{i}, (s^{j})_{j}}}\
        s^{0}
    \end{equation}
    \begin{align}
        \textbf{subject to}\
         & x^0 {}={} x,
        x^{i_+} {}={} A^{i_+}x^i + B^{i_+}u^i,          \\
         & y^i \succcurlyeq_{(\K^i)^*} 0,
        y\iTop b^i {}\leq{} s^i,                        \\
         & E\iTop y^i {}={} \tau^{[i]} + s^{[i]},
        F\iTop y^i {}={} 0,                             \\
         & \ell^{i_+}(x^i,u^i,w^{i_+}) \leq \tau^{i_+},
        \ell^j_N(x^{j}) \leq s^{j},                     \\
         & \Gamma_x^i x^i + \Gamma_u^i u^i \in C^i,
        \Gamma_N^j x^j \in C_N^j,
        \label{eq:constraint_2}
    \end{align}%
\end{subequations}%
for $t\in\N_{[0,N-1]}, i\in\nodes(t), i_+\in\ch(i)$, and $j\in\nodes(N)$, where $(\cdot)^{[i]} = ((\cdot)^{i_+})_{i_+\in\ch(i)}$, $\tau$ and $s$ are slack variables, and $x$ is the initial state.

\subsection{Operator splitting}
\label{sec:operator-splitting}
Define $z = (s^0, z_1, z_2)$, where
$
    z_1 = ((x^{i'})_{i'}, (u^i)_{i}),
$ and $
    z_2 = (y^i,\tau^{[i]},s^{[i]}),
$
for $i'\in \nodes(0, N)$ and $i \in \nodes(0, N-1)$.
We define the sets
\begin{align}
    \mS_1
    :=         &
    \left\{
    z_1
    \left|
    \begin{array}{l}
        x^0 - x = 0,
        \\
        x^{i_+} - A^{i_+} x^{i} - B^{i_+} u^{i} = 0,
        \\
        \text{for } i\in\nodes(0,N-1), i_+\in\ch(i).
    \end{array}
    \right.
    \right\},
    \\
    \mS_2
    :=         & \prod_{i\in\nodes{(0,N-1)}}
    \underbracket[0.5pt]{
        \ker
        \begin{bmatrix}
            E\iTop & -I & -I
            \\
            F\iTop & 0  & 0
        \end{bmatrix}
    }_{\mS_2^i},
    \\ \label{eq:S3-def}
    \mS_3
    :=         &
    \prod_{i\in\nodes{(0,N-1)}}
    \left((\K^i)^* \times \R_+ \times C^i\right)
    \notag
    \\
    {}\times{} &
    \prod_{i\in\nodes{(1,N)}}
    \left( \SOC_{n_x + n_u + 2} + \begin{bmatrix}
        0_{n_x + n_u}\Top & \frac{1}{2} & -\frac{1}{2}
    \end{bmatrix}\Top
    \right)
    \notag
    \\
    {}\times{} &
    \hspace*{-1.2em}\prod_{j\in\nodes{(N)}}
    \left(C^j
    \times
    \hspace*{-0.2em}\left(
        \SOC_{n_x + 2} + \begin{bmatrix}
            0_{n_x}\Top & \frac{1}{2} & -\frac{1}{2}
        \end{bmatrix}\Top
        \right)
    \right).
\end{align}

Next, we define a linear operator $L$ that maps $z$ to
\begin{align}
    \eta :=
    (
     & (
    y^i,
    s^i - b\iTop y^i,
    \Gamma_x^{i} x^i + \Gamma_u^{i} u^i
    )_{i\in\nodes(0,N-1)},
    \notag
    \\
     & (
    ( Q^i )^{\nicefrac{1}{2}} x^{\anc(i)},
    ( R^i )^{\nicefrac{1}{2}} u^{\anc(i)},
    \tfrac{1}{2}\tau^{i},
    \tfrac{1}{2}\tau^{i}
    )_{i\in\nodes(1, N)},
    \notag
    \\
     & (
    \Gamma_N^j x^j,
    (( Q_N^j )^{\nicefrac{1}{2}} x^j,
    \tfrac{1}{2}s^j,
    \tfrac{1}{2}s^j)
    )_{j\in\nodes(N)}
    ).
\end{align}

Remark that $L$ can be made blockdiagonal by permuting its columns, such that each block corresponds to the constraints imposed on a single node. The permutation operation does not change the matrix norm, and the norm of a blockdiagonal matrix equals the largest norm of its constituting blocks. Therefore, $\|L\|$ does not scale with the prediction horizon $N$ and, importantly, neither does the CP step size $\alpha$.

We can now write Problem \eqref{eq:min-s0-v1} as
\begin{subequations}
    \begin{align}
        \minimise_{z}\  & s^0
        \\
        \textbf{subject to}\
                        & z_1 \in \mS_1,
        z_2 \in \mS_2,
        \eta \in \mS_3,
    \end{align}
\end{subequations}
which is equivalent to
\begin{equation}
    \minimise_{z}\ f(z) + g(Lz),
    \label{eq:cp-general-splitting}
\end{equation}
where
\begin{subequations}
    \begin{align}
         & f(z) {}:={}
        s^0
        {}+{}
        \delta_{\mS_1}(z_1)
        {}+{}
        \delta_{\mS_2}(z_2),
        \label{eq:f}
        \\
         & g(\eta) {}:={} \delta_{\mS_3}(\eta).
        \label{eq:g}
    \end{align}
\end{subequations}

\subsection{Proximal operator of $f$}
From Equation \eqref{eq:f}, the separable sum property \cite[Sec.~2.1]{parikh_proximal}, and the fact that $\prox_{\alpha\delta_\mS}(\cdot) = \proj_{\mS}(\cdot)$, we can compute the proximal operator of $f$ by independently computing
$\prox_{\alpha \id}(s^0)$ \cite[Sec. 6.1.1]{parikh_proximal},
$\proj_{\mS_1}(z_1)$, and
$\proj_{\mS_2}(z_2)$.

\subsubsection{Projection on $\mS_1$.}
The projection of $\bar{z}_1$ on $\mS_1$ is the minimizer of the problem
\begin{multline}
    \minimise_{z_1 \in \mS_1}
    \Bigg\{
    \underbracket[0.5pt]{
        \sum_{i\in\nodes(0,N-1)}
        \hspace{-1em}
        \tfrac{1}{2}\left(
        \|x^i-\xb^i\|^2 + \|u^i-\ub^i\|^2
        \right)}_{\text{stage-wise sq. distances}}
    \\
    +
    \underbracket[0.5pt]{
        \sum_{j\in\nodes(N)}
        \tfrac{1}{2}\|x^j-\xb^j\|^2
    }_{\text{terminal stage sq. distances}}
    \Bigg\}.
\end{multline}
This follows the structure of a finite horizon linear-quadratic optimal control (FHOC) problem, which can be solved using the dynamic programming (DP) method. The method for computing $\proj_{\mS_1}(\zb_1)$ comprises the offline Algorithm \ref{alg:s1projection:offline}, and the online Algorithm \ref{alg:s1projection:online}. The offline algorithm is run once at the start and the online algorithm is then run every time $\proj_{\mS_1}(\zb_1)$ is computed.

\begin{algorithm}[htbp!]
    \caption{Projection on $\mS_1$: Offline}\label{alg:s1projection:offline}
    \begin{algorithmic}
        \Require the system matrices $A{\in}\R^{n_x\times n_x}$, $B{\in}\R^{n_x\times n_u}$, and the prediction horizon $N$

        \Ensure Matrices $(K^i)_i$, $(P^i)_i$, $(\widetilde{R}^i)_i$ and $(\bar{A}^i)_i$

        \ForAll{$i\in\nodes(N)$ \textbf{in parallel}}
        \State $P^i \gets I_{n_x}$
        \EndFor

        \For{$t=0,1,\ldots, N-1$}
        \ForAll{$i\in\nodes(N-(t+1)) \textbf{ in parallel}$}
        \State$
            \widetilde{P}^{i_+} \gets
            B^{i_+\intercal} P^{i_+}
        $, $i_+{\in}\ch(i)$
        \State$
            \widetilde{R}^i \gets $\parbox[t]{170pt}{$
                I_{n_u}
                + \sum_{i_+\in\ch(i)} \widetilde{P}^{i_+} B^{i_+}
            $
            and determine its Cholesky decomposition \strut}
        \State$
            K^i \gets
            -(\widetilde{R}^{i})^{-1}
            \sum_{i_+\in\ch(i)} \widetilde{P}^{i_+} A^{i_+}
        $
        \State$
            \bar{A}^{i_+} \gets
            A^{i_+} + B^{i_+}K^i
        $, $i_+{\in}\ch(i)$
        \State$
            P^i \gets
            I_{n_x} {+} K^{i\intercal} K^i {+} \sum_{i_+\in\ch(i)} \bar{A}^{i_+\intercal} P^{i_+} \bar{A}^{i_+}
        $
        \EndFor

        \EndFor

    \end{algorithmic}
\end{algorithm}

\begin{algorithm}[htbp!]
    \caption{Projection on $\mS_1$: Online}\label{alg:s1projection:online}
    \begin{algorithmic}
        \Require Matrices computed by Algorithm \ref{alg:s1projection:offline}, the vector $\zb_1$ as defined in Sec.~\ref{sec:operator-splitting}, and the initial state $x^0$

        \Ensure Projection on $\mS_1$ at $\zb_1$

        \ForAll{$i\in\nodes(N)$ \textbf{in parallel}}
        \State $q^i_0 \gets -\xb^i$
        \EndFor

        \For{$t=0,1,\ldots, N-1$}
        \ForAll{$i\in\nodes(N-(t+1)) \textbf{ in parallel}$}
        \State$
            d_{t+1}^i \gets $\parbox[t]{170pt}{$
                (\widetilde{R}^{i})^{-1}
                \left(
                \ub^i - \sum_{i_+\in\ch(i)} B^{i_+\intercal} q^{i_+}_{t}
                \right)
            $
            using the Cholesky decomposition of $\widetilde{R}^{i}$\strut}
        \State$
            q_{t+1}^i {\gets} $\parbox[t]{0pt}{$
                K^{i\intercal} \left( d^i {-} \ub^i \right)
                {-} \xb^i
                + \sum_{i_+\in\ch(i)} \bar{A}^{i_+\intercal} \left( P^{i_+} B^{i_+} d^i + q^{i_+} \right)
            $\strut}
        \EndFor
        \EndFor

        \State $x^0 \gets x$
        \For{$t=0,1,\ldots, N-1$}
        \ForAll{$i\in\nodes(t), i_+\in\ch(i)$ \textbf{in parallel}}
        \State $u^i \gets K^i x^i + d^i$
        \State $x^{i_+} \gets A^{i_+}x^i + B^{i_+}u^i$
        \EndFor
        \EndFor

        \State \Return $z_1 = \proj_{\mS_1}(\zb_1)$.

    \end{algorithmic}
\end{algorithm}

\subsubsection{Projection on $\mS_2$.}
Since $\mS_2 = \prod_{i\in\nodes(0,N-1)} \mS_2^i$, the projection on $\mS_2$ is the projection on the kernel of a  matrix of typically small dimensions, i.e., the projection on a linear space. This can either be solved very efficiently and accurately using a numerical method, or, by precomputing a pseudoinverse. Note that the projections for each scenario tree node can be computed in parallel.

\subsection{Proximal operator of $g^*$}
We start with the proximal operator of $g$ with parameter $\alpha$, where
$
    \prox_{\alpha g}(\eta)
    {}={}
    \proj_{\mS_3}(\eta).
$
Then, by the extended Moreau decomposition \cite[Thm.6.45]{beck_optimization}
\begin{equation}\label{eq:prox:g}
    \prox_{\alpha g^*}(\eta)
    {}={} \eta - \proj_{\mS_3} \left( \tfrac{\eta}{\alpha} \right).
\end{equation}
The set $\mS_3$ is the Cartesian product of low-dimensional (translated) convex cones (cf. \eqref{eq:S3-def}),
projections on which can be carried out independently.

\section{\SPOCK: SuperMann and Chambolle-Pock}\label{sec:spock}
To accelerate the CP method we combine it with SuperMann \citep{themelis_supermann} and call the resulting ensemble \SPOCK{} (see Algorithm \ref{alg:spock}). SuperMann is an algorithmic framework that augments the classical Krasnosel'ski\v{\i}-Mann iteration of an FNE operator with a line search to yield a method that enjoys the same global convergence properties, and under certain weak assumptions converges superlinearly.
The FNE operator for \SPOCK{} is $T$ defined by Equation \eqref{eq:chambolle-pock-operator}.
We combine SuperMann with Anderson's acceleration (AA) method \cite[]{anderson_acceleration_original} to compute fast quasi-Newtonian directions.
The most computationally demanding parts of the algorithm are the operators $L$ and $L^*$, which are computed in $T$.

Let us first define $c_{z}^{(k)} := z^{(k)} - z^{(k+1)}$, and likewise we define $c_{\eta}^{(k)}$. We denote a primal update direction by $d_{z}$, and a dual update direction by $d_{\eta}$. Now we can define $v^{(k)} {:=} (z^{(k)}, \eta^{(k)})$, $c^{(k)} {:=} (c_{z}^{(k)}, c_{\eta}^{(k)})$, and $d^{(k)} {:=} (d_{z}^{(k)}, d_{\eta}^{(k)})$.

\subsection{SuperMann}
SuperMann is a Newton-type algorithm that finds a fixed-point $v^{\star}\in\fix T$ by finding a zero of $R$. The framework involves two extragradient-type updates of the general form
\begin{subequations}
    \begin{equation}\label{eq:K1}
        v_{\textit{K1}}^{(k+1)} = v^{(k)} + \tau_{(k)} d^{(k)},
    \end{equation}
    \begin{equation}\label{eq:K2}
        v_{\textit{K2}}^{(k+1)} = v^{(k)} - \bar{\tau}_{(k)} R(v_{\textit{K1}}^{(k+1)}),
    \end{equation}
\end{subequations}
where $d^{(k)}$ are fast, e.g., quasi-Newtonian, directions and the step sizes $\tau_{(k)}$ and $\bar{\tau}_{(k)}$ are chosen so that SuperMann enjoys the same global convergence properties as the classical Krasnosel'ski\v{\i}-Mann scheme. Note that when the directions satisfy a Dennis-Mor\'{e} condition, SuperMann converges superlinearly under mild assumptions that include metric subregularity of the residual.

At each step we perform a backtracking line search on $\tau_{(k)}$ until we either trigger an ``educated'' \textit{K1} update, (fast convergence, cf. \eqref{eq:K1}) or a ``safeguard'' \textit{K2} update, (global convergence, cf. \eqref{eq:K2}) as shown in Algorithm \ref{alg:spock}. The \textit{K2} update can be interpreted as a projection of $v^{(k)}$ on a hyperplane generated by $v_{\textit{K1}}^{(k+1)}$ (which separates $\fix T$ from $v^{(k)}$), therefore guaranteeing that every iterate $v^{(k+1)}$ moves closer to $\fix T$. SuperMann also checks for a sufficient decrease of the norm of the residual $\|R(v^{(k)})\|$ at each step, which triggers a ``blind'' \textit{K0} update of the form $v^{(k+1)}_{\textit{K0}} = v^{(k)} + d^{(k)}$. This does not require line search iterations.

\subsection{Anderson's acceleration}\label{sec:aa}
Quasi-Newtonian directions can be computed according to the general rule
\begin{equation}
    d^{(k)} = -B_{(k)}^{-1}c^{(k)} = -H_{(k)}c^{(k)},
\end{equation}
where the invertible linear operators $H_{(k)}$ are updated according to certain low-rank updates so as to satisfy certain secant conditions, starting from an initial operator $H_{(0)}$. AA enforces a multi-secant condition \cite[]{fang_multisecant, walker_anderson_acceleration}. For each iteration $k$, we update two buffers of length $m$: one for the differences of $z$ and $\eta$ between iterations,
\begin{equation}
    M_{(k)}^{\rm{P}} =
    \begin{bmatrix}
        v^{(k)}-v^{(k-1)} \cdots v^{(k-m+2)}-v^{(k-m+1)}
    \end{bmatrix},
\end{equation}
and one for the differences of $c_{z}$ and $c_{\eta}$ between iterations,
\begin{equation}
    M_{(k)}^{\rm{R}} =
    \begin{bmatrix}
        c^{(k)}-c^{(k-1)} \cdots c^{(k-m+2)}-c^{(k-m+1)}
    \end{bmatrix}.
\end{equation}
Directions are computed by
\begin{equation}
    d^{(k)} = -c^{(k)} - (M_{(k)}^{\rm{P}}-M_{(k)}^{\rm{R}})\gamma^{(k)},
\end{equation}
where $\gamma^{(k)}$ is the least squares solution to the linear system $M_{(k)}^{R}\gamma^{(k)} = c^{(k)}$, that is, $\gamma^{(k)}$ solves
\begin{equation}
    \minimise_{\gamma^{(k)}} \|M_{(k)}^{R}\gamma^{(k)} - c^{(k)}\|^2,
\end{equation}
which can be solved by QR factorization of the matrix $M_{(k)}^{R}$ (which may be updated at every iteration \cite[]{walker_anderson_acceleration}).
In practice AA performs well with short memory lengths between 3 and 10.

\begin{algorithm}
    \caption{\SPOCK{} algorithm for RAOCPs}\label{alg:spock}
    \begin{algorithmic}[1]
        \Require
        problem data,
        $z^{(0)}$ and $\eta^{(0)}$,
        tolerances $\epsilon_{\text{abs}} > 0$ and $\epsilon_{\text{rel}} > 0$,
        $\alpha$ such that $0 < \alpha \|L\| < 1$,
        $\N \ni m > 0$,
        $c_0, c_1, c_2 \in [0,1)$,
        $\beta,\sigma \in (0,1)$, and
        $\lambda \in (0,2)$.

        \Ensure approximate solution of RAOCP

        \State Invoke Algorithm \ref{alg:s1projection:offline}

        \State $r^{(0)} \gets v^{(0)} - T(v^{(0)})$

        \State $\zeta_{(0)} {\gets} ( (r_{z}^{(0)})\Top (r_{z}^{(0)} - \alpha L^{\ast}r_{\eta}^{(0)}) {+} (r_{\eta}^{(0)})\Top (r_{\eta}^{(0)} - \alpha Lr_{z}^{(0)}) )^{\nicefrac{1}{2}}$

        \State $\omega_{\text{safe}} \gets \zeta_{(0)}$, $k \gets 0$

        \If{termination criteria (Sec.~\ref{sec:cp-method}) are satisfied}
        \State \textbf{return} $v^{(k)}$\label{line:check-term}
        \EndIf

        \State Use AA to compute an update direction, $d^{(k)}$ \hfill(Sec.\ref{sec:aa})

        \State $r^{(k)} \gets v^{(k)} - T(v^{(k)})$

        \State $\omega_{(k)} {\gets} ((r_{z}^{(k)})\Top (r_{z}^{(k)} {-} \alpha L^{\ast}r_{\eta}^{(k)}) {+} (r_{\eta}^{(k)})\Top (r_{\eta}^{(k)} {-} \alpha Lr_{z}^{(k)}))^{\nicefrac{1}{2}}$

        \If{$\omega_{(k)} \leq c_0 \zeta_{(k)}$}
        \State $v^{(k+1)} {\gets} v^{(k)} {+} d^{(k)}$, $\zeta_{(k+1)} {\gets} \omega_{(k)}$, \textbf{goto} line \ref{line:last-step} \hfill{\color{gray}(\textit{K0})}
        \EndIf

        \State $\zeta_{(k+1)} \gets \zeta_{(k)}$, $\tau \gets 1$

        \State $\tilde{v}^{(k)} \gets v^{(k)} + \tau d^{(k)}$, $\tilde{r}^{(k)} \gets \tilde{v}^{(k)} - T(\tilde{v}^{(k)})$ \label{line:try-new-tau}

        \State $\tilde{\omega}_{(k)} {\gets} ((\tilde{r}_{z}^{(k)})\Top (\tilde{r}_{z}^{(k)} {-} \alpha L^{\ast}\tilde{r}_{\eta}^{(k)}) {+} (\tilde{r}_{\eta}^{(k)})\Top (\tilde{r}_{\eta}^{(k)} {-} \alpha L\tilde{r}_{z}^{(k)}))^{\nicefrac{1}{2}}$

        \If{$\omega_{(k)} \leq \omega_{\text{safe}}$ \textbf{and} $\tilde{\omega}_{(k)} \leq c_1 \omega_{(k)}$}
        \State $v^{(k+1)} {\gets} \tilde{v}^{(k)}$, $\omega_{\text{safe}} {\gets} \tilde{\omega}_{(k)} {+} c_2^k$, \textbf{goto} line \ref{line:last-step} \hfill{\color{gray}(\textit{K1})}
        \EndIf

        \State $\rho_{(k)} \gets \tilde{\omega}_{(k)}^2 - 2\alpha (\tilde{r}^{(k)})\Top M ( \tilde{v}^{(k)} - v^{(k)} )$

        \If{$\rho_{(k)} \geq \sigma \tilde{\omega}_{(k)} \omega_{(k)}$}
        $v^{(k+1)} {\gets} v^{(k)} {-} \dfrac{\lambda \rho_{(k)}}{\tilde{\omega}_{(k)}^2} \tilde{r}^{(k)}$ \hfill{\color{gray}(\textit{K2})}
        \EndIf
        \State \textbf{else} $\tau \gets \beta \tau$, \textbf{goto} line \ref{line:try-new-tau}

        \State $k \gets k + 1$, \textbf{goto} line \ref{line:check-term} \label{line:last-step}
    \end{algorithmic}
\end{algorithm}

\section{Simulation results}

Consider a linear system $x_{t+1} = A(w_t) x_t + B(w_t) u_t$ where $A(w)$ is the tridiagonal matrix defined by
\begin{subequations}\label{eq:dynamics-example}
    \begin{align}
         & A_{j, j}(w) = 1 + \frac{w - 1}{d} \left(1 + \frac{j - 1}{n_x}\right), &  & j \in \N_{[1, n_x]}, \label{eq:dynamics-example-diag}    \\
         & A_{j, j - 1}(w) = A_{j - 1, j}(w) = 0.01,                             &  & j \in \N_{[2, n_x]}, \label{eq:dynamics-example-supdiag}
    \end{align}
\end{subequations}
and $B(w) = I_{n_u}$ for $w\in\N_{[1, d]}$.
This system is inspired by \cite[]{schuurmans_example,dean_sample-complexity-lqr} and models the temperature of $n_x$ servers in a data center. The state $x_t \in \R^{n_x}$ corresponds to the deviations from some nominal temperature, and the input $u_t \in \R^{n_u}$ models the amount of heating ($u_t \geq 0$) or cooling ($u_t < 0$) that is applied to each server.
The diagonal elements \eqref{eq:dynamics-example-diag},
and the sub- and super-diagonal elements \eqref{eq:dynamics-example-supdiag} of $A(w)$ model the heat generated by each individual machine, and the heat transfer between adjacent machines, respectively.
The disturbance $w$ represents the load on the system; if $w = 1$, the servers are idle and no additional heat is generated. If $w = d$, all servers are running at full capacity and a maximum amount of heat is generated.
Note that server $j$ produces more heat than server $j - 1$ under the same load.

We use the quadratic stage and terminal costs from Sec.~\ref{sec:scenario-trees} where $Q^{i_+}=Q^j_N=I_{n_x}$ and $ R^{i_+}=10I_{n_u}$ for $i\in\nodes(0,N-1), i_+\in\ch(i)$ and $j\in\nodes(N)$. We also impose the constraints $\|x^i\|_{\infty} \leq 1$ for $i\in\nodes(0,N)$ and $\|u^i\|_{\infty} \leq 1.5$ for $i\in\nodes(0,N-1)$. Each node has an associated $\AVaR_{a}$ with $a = 0.95$. This section compares \SPOCK\footnote{A sequential implementation is available at \href{https://github.com/kul-optec/spock.jl}{https://github.com/kul-optec/spock.jl}.} to existing solvers when solving the RAOCP (\ref{eq:original-problem}) of horizon $N$ for the above example. Remark that the number of scenarios in this problem, $d^N$, grows exponentially with the prediction horizon, $N$.

In what follows, we fix $d = 2$ and consider reference probabilities
$\pi^i = ( 0.3, 0.7 )\Top$, for $i\in\nodes(0, N-1)$.
The results for \SPOCK{} are obtained using the SuperMann parameters from \cite[Sec.VI.D]{themelis_supermann}. Quasi-Newton directions are generated by AA with $m = 3$. The other solvers all solve the tractable reformulation (\ref{eq:min-s0-v1}) through the \textsc{JuMP} interface \citep{jump}.

\subsection{SuperMann acceleration}
To illustrate the benefit of applying SuperMann, we compare the performance of the vanilla CP algorithm with the \SPOCK{} solver. Consider system \eqref{eq:dynamics-example} with $N = 7$, $n_x = 5$, and $x^0 = 0.1 \cdot 1_{n_x}$.

\begin{figure}
    \centering
    \resizebox{\columnwidth}{!}{
        \input{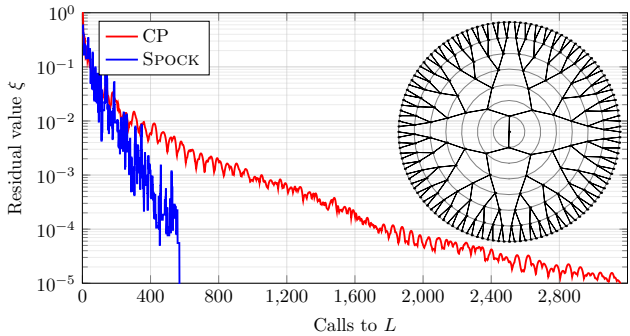}
    }
    \caption{Residual value $\xi$ against the number of calls of the linear operator for the vanilla Chambolle-Pock algorithm (CP) and \SPOCK, for $N = 7, n_x = 5$. A visualization of the scenario tree is shown, in which the root node is at the center and each concentric circle corresponds to a time stage.}
    \label{fig:cp-vs-cpsp}
\end{figure}

Figure \ref{fig:cp-vs-cpsp} shows the residual value $\xi$ from Sec.~\ref{sec:cp-method} against the number of calls
to the operator
$L$ for both CP and \SPOCK, with tolerances $\epsilon_{\text{abs}} = \epsilon_{\text{rel}} = 10^{-5}$.
It demonstrates that SuperMann requires considerably fewer calls compared to vanilla CP,
as the former converges with 571 calls, and the latter with 3159.
Furthermore, Fig.~\ref{fig:cp-vs-cpsp} suggests that this difference becomes
more pronounced as the tolerances decrease.

\subsection{Open-loop simulation}\label{sec:open-loop}
We now evaluate \SPOCK's performance by comparing it to both open-source solvers (\IPOPT{} \cite[]{ipopt}, \SEDUMI{} \cite[]{sedumi}, \COSMO{} \cite[]{cosmo}) and proprietary alternatives (\MOSEK{} \cite[]{andersen_mosek}, \GUROBI{} \cite[]{gurobi_manual}). \COSMO{} and \SEDUMI{} are conic solvers, while \IPOPT, \MOSEK{} and \GUROBI{} are interior-point solvers.

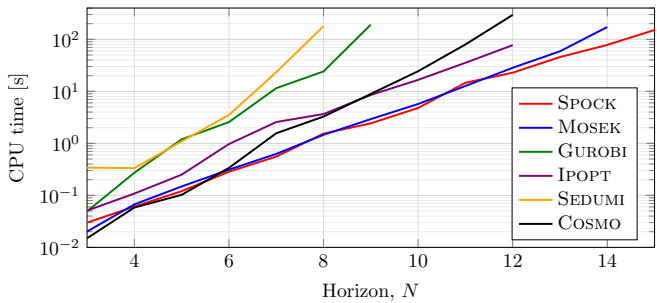
\begin{figure}
    \centering
    \resizebox{\columnwidth}{!}{
        \begin{tikzpicture}
    \begin{axis}[
            width=4in, height=1.7in,
            at={(1.011in,0.642in)},
            legend cell align={left},
            legend pos=south east,
            scale only axis,
            minor grid style={thin,draw opacity=0.3},
            major grid style={thin,draw opacity=0.5},
            grid=both,
            xmin=3,
            xmax=15,
            ymin=0.01,
            ymax=400,
            ymode=log,
            xlabel = {Horizon, $N$},
            ylabel = {CPU time [s]}]
            \addplot[
                color={rgb,1:red,1.0;green,0.0;blue,0.0},
                name path={7af44605-0608-4499-a7b4-aefbe091aa5c}, draw opacity={1.0}, line width={1}, solid
            ]
            table[row sep={\\}]
                {
                    \\
                    3.0  0.02994837  \\
                    4.0  0.060908716  \\
                    5.0  0.120554163  \\
                    6.0  0.28526187  \\
                    7.0  0.557456879  \\
                    8.0  1.540880604  \\
                    9.0  2.428256747  \\
                    10.0  4.800994611  \\
                    11.0  14.616935775  \\
                    12.0  22.880396134  \\
                    13.0  45.661248555  \\
                    14.0  78.266806084  \\
                    15.0  151.511085064  \\
                }
                ;
            \addlegendentry {\SPOCK}
            \addplot[color={rgb,1:red,0.0;green,0.0;blue,1.0}, name path={c0742d21-27cb-4820-be90-577c68efae52}, draw opacity={1.0}, line width={1}, solid]
                table[row sep={\\}]
                {
                    \\
                    3.0  0.02019672  \\
                    4.0  0.067015278  \\
                    5.0  0.148361053  \\
                    6.0  0.311760679  \\
                    7.0  0.627032835  \\
                    8.0  1.460537915  \\
                    9.0  2.929311131  \\
                    10.0  5.712879452  \\
                    11.0  12.666745256  \\
                    12.0  28.423045533  \\
                    13.0  58.958225147  \\
                    14.0  172.560179654  \\
                }
                ;
            \addlegendentry {\MOSEK}
            \addplot[color={rgb,1:red,0.0;green,0.502;blue,0.0}, name path={a128064e-38bf-4a1d-8a0c-15bb103e4380}, draw opacity={1.0}, line width={1}, solid]
                table[row sep={\\}]
                {
                    \\
                    3.0  0.049411809  \\
                    4.0  0.272411627  \\
                    5.0  1.196854314  \\
                    6.0  2.564091705  \\
                    7.0  11.504945539  \\
                    8.0  24.149949576  \\
                    9.0  191.923651687  \\
                }
                ;
            \addlegendentry {\GUROBI}
            \addplot[color={rgb,1:red,0.502;green,0.0;blue,0.502}, name path={277f12e5-720c-45c5-9d66-e311f08e5287}, draw opacity={1.0}, line width={1}, solid]
                table[row sep={\\}]
                {
                    \\
                    3.0  0.050937591  \\
                    4.0  0.108340217  \\
                    5.0  0.250095871  \\
                    6.0  0.970372604  \\
                    7.0  2.572118021  \\
                    8.0  3.6617332  \\
                    9.0  8.532633751  \\
                    10.0  16.620193589  \\
                    11.0  35.418323034  \\
                    12.0  77.359181342  \\
                }
                ;
            \addlegendentry {\IPOPT}
            \addplot[color={rgb,1:red,1.0;green,0.6471;blue,0.0}, name path={ce69c858-dff9-43d6-9e37-3dd2d7b9c80c}, draw opacity={1.0}, line width={1}, solid]
                table[row sep={\\}]
                {
                    \\
                    3.0  0.341806092  \\
                    4.0  0.333417755  \\
                    5.0  1.092254957  \\
                    6.0  3.523998985  \\
                    7.0  23.339722702  \\
                    8.0  179.297366053  \\
                }
                ;
            \addlegendentry {\SEDUMI}
            \addplot[color={rgb,1:red,0.0;green,0.0;blue,0.0}, name path={9827b513-4f3d-4fef-bded-88465a488968}, draw opacity={1.0}, line width={1}, solid]
                table[row sep={\\}]
                {
                    \\
                    3.0  0.014994549  \\
                    4.0  0.058373584  \\
                    5.0  0.102449652  \\
                    6.0  0.337407834  \\
                    7.0  1.564327529  \\
                    8.0  3.273527758  \\
                    9.0  8.923604297  \\
                    10.0  24.39213591  \\
                    11.0  79.657784888  \\
                    12.0  294.944967293  \\
                }
            ;
            \addlegendentry {\COSMO}
    \end{axis}
\end{tikzpicture}
    }
    \caption{Execution time against the prediction horizon $N$, for $n_x=n_u=50$. The number of scenarios is $2^N$.}
    \label{fig:scaling}
\end{figure}

Consider system \eqref{eq:dynamics-example} with $n_x = 50$, $x^0 = 0.1 \cdot 1_{n_x}$, and $\epsilon_{\text{abs}} = \epsilon_{\text{rel}} = 10^{-3}$.
Figure \ref{fig:scaling} shows the execution time of the aforementioned solvers against the prediction horizon $N$.
The number of iterations of \SPOCK{} in this simulation is shown in Table \ref{tab:iterations}.
Once a solver reaches the run time threshold of \unit[150]{s}, the horizon $N$ is no longer increased for that solver.
The performances of \SEDUMI{} and \GUROBI{} appear to rapidly degrade as the prediction horizon increases.
\IPOPT{} and \COSMO{} perform decently for moderate problem sizes.
However, for large values of $N$, the execution time of \COSMO{} increases at a considerably higher rate than that of \SPOCK, and \IPOPT{} even runs out of memory for $N = 12$.
\SPOCK{} performs competitively with \MOSEK{} for smaller problem sizes, and even outperforms it by approx. a factor 2   for $N=14$.
Further experiments have indicated that eventually \MOSEK{} runs out of memory when $N = 16$, whereas \SPOCK{} does not.
Altogether, these results are highly promising, considering that the computationally demanding steps of the \SPOCK{} solver can be massively parallelized.
\begin{table}[ht!]
    \centering
    \caption{\SPOCK{} iterations corresponding to Fig.~\ref{fig:scaling}.}
    \label{tab:iterations}
    \begin{tabular}{@{}cccccc@{}}
        \toprule
        Horizon, $N$    & 3   & 6   & 9   & 12  & 15  \\ \midrule
        Num. iterations & 105 & 128 & 110 & 109 & 105 \\ \bottomrule
    \end{tabular}
\end{table}

\subsection{Closed-loop simulation}
Most of the solvers we compared \SPOCK{} to use interior-point methods and are, therefore, difficult to warm-start.
However, within the context of MPC, a reasonable initial guess of the solution is typically available from the previous time step.
We demonstrate that this can be leveraged to a far greater extent by the \SPOCK{} solver than by the compared alternatives, by carrying out
an MPC simulation of $20$ steps, in which we iteratively solve the RAOCP for successive initial states.
Warm-starts are performed by using the primal and dual solutions at one MPC step to initialize the primal and dual vectors in the next MPC step.
We consider system \eqref{eq:dynamics-example} with $N = 10$, $n_x = 20$, $\epsilon_{\text{abs}} = \epsilon_{\text{rel}} = 10^{-3}$ and the initial state in the first MPC step $x^0 = 0.1 \cdot 1_{n_x}$.

\begin{figure}
    \centering
    \resizebox{\columnwidth}{!}{
        \input{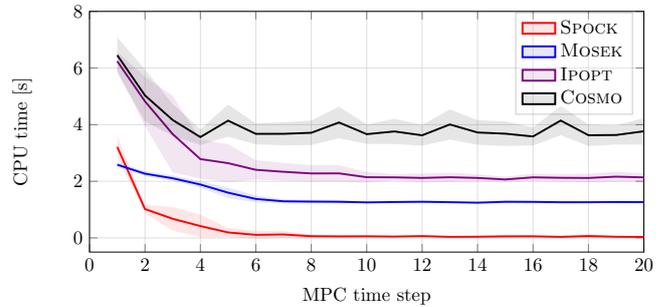}
    }
    \caption{Average execution time for successive MPC time steps over 15 random realizations with $N = 10$ and $n_x = n_u = 20$. \SPOCK{} and \COSMO{} exploit warm-starts, the interior-point solvers do not.}
    \label{fig:mpc}
\end{figure}

Figure \ref{fig:mpc} compares the average execution time of the best performing solvers mentioned in Sec.~\ref{sec:open-loop}, for successive MPC time steps. \SPOCK{} and \COSMO{} exploit warm-starts, the interior-point solvers do not. In the first MPC step, \SPOCK{} performs better than \IPOPT{} and \COSMO, and slightly slower than \MOSEK. However, warm-starts greatly reduce the required number of iterations for \SPOCK, such that it outperforms all the other solvers in closed loop.%

\section{Conclusions and future work}
We proposed a splitting scheme for general risk-averse problems tailored to proximal methods by exploiting the problem structure.
The splitting is used in our proposed algorithm \SPOCK{} which is amenable to warm-starting and massive parallelization.
\SPOCK{} takes advantage of the SuperMann scheme combined with fast directions from Anderson's acceleration method to enhance the convergence speed. We implemented \SPOCK{} in \texttt{Julia} as an open-source solver, which is available on GitHub at \href{https://github.com/kul-optec/spock.jl}{https://github.com/kul-optec/spock.jl}.
Our simulations demonstrate: the benefit of \SPOCK{} over the vanilla CP method; 
that the solver, even without parallelization, outperforms proprietary solvers such as \MOSEK{} and \GUROBI for large problem sizes;
and, that when warm-starting is used in the context of MPC, \SPOCK{} outperforms the compared solvers by a considerable margin.
In the future, we plan to harness the parallelization capabilities of GPUs by implementing \SPOCK{} in CUDA-C++.

\bibliography{root.bib}

\end{document}